\documentclass[11pt]{amsart}
\usepackage{amsmath}
\usepackage{amsfonts,amssymb}
\usepackage{amsthm}

\usepackage[usenames]{xcolor}
\usepackage{color}
\usepackage{caption}
\usepackage{graphics}
\usepackage{mathtools}
\usepackage{diagbox}

\renewcommand{\epsilon}{\varepsilon}
\newcommand{\E}{\mathbb{E}} 
\newcommand{\N}{\mathbb{N}}
\newcommand{\C}{\mathbb{C}}
\newcommand{\R}{\mathbb{R}}

\newcommand{\ntr}{\mathbf{tr}}
\newcommand{\diag}{\textrm{diag}}

\newtheorem{theorem}{Theorem}
\newtheorem{lemma}[theorem]{Lemma}
\newtheorem{corollary}[theorem]{Corollary}
\newtheorem{proposition}[theorem]{Proposition}
\newtheorem{problem}[theorem]{Problem}

\theoremstyle{definition}

\newtheorem{example}[theorem]{Example}
\newtheorem{remark}[theorem]{Remark}

\newcommand{\be}{\begin{eqnarray}}
\newcommand{\ee}{\end{eqnarray}}
\newcommand{\ben}{\begin{enumerate}}
\newcommand{\een}{\end{enumerate}}

\newcommand{\ba}{\begin{array}}
\newcommand{\ea}{\end{array}}

\newcommand{\mc}{\mathcal}

\newcommand{\tr}{\mathrm{tr}}
\newcommand{\Her}{\mathrm{Her}}
\newcommand{\spec}{\mathrm{sp}}

\newcommand{\eps}{\varepsilon}

\newcommand\xqed[1]{%
  \leavevmode\unskip\penalty9999 \hbox{}\nobreak\hfill
  \quad\hbox{#1}}
\newcommand\demo{\xqed{$\triangle$}}

\begin{document}

\title{{Optimal bounds on the positivity of a matrix from a few moments}}

\author{Gemma De las Cuevas}
\address{Institute of Theoretical Physics, Technikerstr.\ 21a, A6020 Innsbruck, Austria}
\email{gemma.delascuevas@uibk.ac.at}

\author{Tobias Fritz}
\address{Max Planck Institute for Mathematics in the Sciences,  Inselstr.\ 22, D04103 Leipzig, Germany}
\email{tobias.fritz@mis.mpg.de}

\author{Tim Netzer}
\address{Department of Mathematics, Technikerstr.\ 13, A6020 Innsbruck, Austria}
\email{tim.netzer@uibk.ac.at}


\maketitle
\begin{abstract}
In many contexts one encounters Hermitian operators $M$ on a Hilbert space whose dimension is so large that it is impossible to write down all matrix entries in an orthonormal basis. How does one determine whether such $M$ is positive semidefinite? 
Here we approach this problem by deriving asymptotically optimal bounds to the distance to the positive semidefinite cone in Schatten $p$-norm for all integer $p\in[1,\infty)$, assuming that we know the moments $\ntr(M^k)$ up to a certain order $k=1,\ldots, m$. 
We then provide three methods to compute these bounds and relaxations thereof: 
the sos polynomial method (a semidefinite program), 
the Handelman method (a linear program relaxation), and the Chebyshev method (a relaxation not involving any optimization). 
We investigate the analytical and numerical performance of these methods and 
present a number of example computations, partly motivated by applications to tensor networks and to the theory of free spectrahedra.
\end{abstract}

\section{Introduction}

\subsection{The setting}

Positive semidefinite matrices and operators come up in a large number of contexts throughout mathematics, physics and computer science, including the theory of operator algebras, computations with quantum states, or semidefinite relaxations of optimization problems. A problem that arises frequently in many of these contexts is determining whether a given matrix $M$ or operator is positive semidefinite (psd). 
In certain applications---we describe two of them in Sections \ref{ssec:appl1} and \ref{ssec:appl2}---$M$ is of such a large dimension that it is impossible to express $M$ in an orthonormal basis and store the result in a computer, let alone diagonalize it or compute determinants of principal minors. Instead, we assume that one can compute a few of the normalized moments  $\ntr(M^k)$, where we define the normalized trace as 
$$
  \ntr(M) := s^{-1}\,\tr(M) ,
$$
where $s$ is the size of $M$. In this paper, we answer the following questions:
\begin{enumerate}
\item[(i)] Given the first $m$ normalized moments $\ntr(M^k)$ for $k=1,\ldots,m$ of a Hermitian operator $M$ with $\|M\|_\infty \leq 1$, can one show that $M$ is not psd?
\item[(ii)] Given these moments and a $p\in [1,\infty)$, can one optimally bound the distance of $M$ to the psd cone from above and below in Schatten $p$-norm? 
\end{enumerate}

Since both the moments and the positive semidefiniteness of a Hermitian operator $M$ are characterized by the distribution of eigenvalues---or more generally by the corresponding Borel measure on the spectrum of $M$---we are secretly concerned with a version of the truncated Hausdorff moment problem. In these terms, the above two questions become:

\begin{enumerate}
\item[(i)] Given only the moments
\[
	\E_\mu[x^k] := \int_{-1}^{+1} x^k \, d\mu
\]
of a compactly supported Borel measure $\mu$ on $[-1,1]$ for $k=1,\ldots,m$, can one show that $\mu$ is not supported on $[0,1]$?
\item[(ii)] Given these moments and a $p\in [1,\infty)$, 
can one optimally bound the $p$-Wasserstein distance \cite{Vi09} between $\mu$ and the set of probability measures supported on $[0,1]$ from above and from below?
\end{enumerate}

Given this connection with the moment problem, it should not come as a surprise that our methods also work in certain infinite-dimensional situations. While we focus on the matrix case in the main text, we sketch the extension to testing positivity of a Hermitian element in a von Neumann algebra equipped with a finite faithful trace in Appendix~\ref{vNa}.


Let us now motivate the assumption that we have access to a few moments of $M$. In the applications we have in mind, 
the space where $M$ lives has a natural tensor product structure, with respect to which $M$ can be expressed as 
 \be \label{eq:Cl}
M=\sum_{j=1}^{r} A^{[1]}_j\otimes \cdots \otimes A^{[n]}_j , 
\ee 
 where each $A^{[i]}_j$ is Hermitian and of reasonably small dimension\footnote{If $r=1$, then there is a simple criterion to determine whether $M$ is psd: $M$ is psd if and only if each $A^{[i]}$ is either psd or negative semidefinite, and the number of negative semidefinite matrices is even. But we are not aware of any such simple criterion for $r>1$.}. Note that every $M$ on a tensor product space can be written this way, for large enough $r$. 
In our example applications below, $r$ is taken to be fixed and typically small, or not scaling with $n$.   Thus, the naturally available operations are those that can be directly performed in terms of the local matrices $A_j^{[i]}$. 
This includes taking powers of $M$ and taking the trace, which gives us access to the moments of $M$:
\begin{equation}\label{eq:tr}
\begin{split}
\ntr\left(M^k\right) 
& = \sum_{j_1,\ldots,j_k=1}^{r} \ntr\left(A_{j_1}^{[1]}\cdots A_{j_k}^{[1]}\right)\cdots\ntr\left(A_{j_1}^{[n]}\cdots A_{j_k}^{[n]}\right),
\end{split}
\end{equation}
for $k\in\N$ much smaller than the size of the matrix $M$, namely $s\times s$.  

\subsection{Example application: tensor networks}
\label{ssec:appl1}

Our first example application concerns quantum states of a many-body system, which  are modelled by psd matrices on a Hilbert space $\mc{H} = \mc{H}_1\otimes \mc{H}_2\otimes \ldots \otimes \mc{H}_n$, where  $\mc{H}_i$ is the Hilbert space associated to subsystem $i$. Typically, all $\mc{H}_i$ are of the same finite dimension. Since the dimension of $\mc{H}$ grows exponentially with $n$, physicists have attempted to develop a  scalable description of quantum many-body systems, that is, one which grows only polynomially in $n$. This is the objective of the program of \emph{tensor networks} \cite{Fa92,Vi03,Pe07,Or14b}. 
While this program has been very successful for pure states (i.e.~psd matrices of rank 1), its success for mixed states has been more limited. 
One of the reasons for that is the \emph{positivity problem}, which is the following. 
In the tensor network paradigm, it is natural to use a few matrices for each local Hilbert space $\mc{H}_i$. For example, the state of the system in one spatial dimension (with periodic boundary conditions) is described by 
\be
M=\sum_{j_1, \ldots , j_{n}=1}^{r}
A_{j_1,j_2} \otimes A_{j_2,j_3} \otimes \cdots \otimes 
A_{j_{n},j_1},
\label{eq:Mtn}
\ee
where each  $A_{j_{l},j_{l+1}}$ is a Hermitian matrix in $\mc{H}_l$ 
\cite{Ve04d,De13c}. 
(In general, the local matrices $A_{j_l,j_{l+1}}$ may also depend on the site $l$, but we do not consider this case for notational simplicity.) 
Now, while $M$ must be psd to describe a quantum state, each of the local matrices $A_{j_{l},j_{l+1}}$ need not be psd. 
While there is a way of imposing positivity in the local matrices (resulting in the `local purification form'), this generally comes at the price of a very large increase in the number of matrices, thus making the representation very inefficient  \cite{De13c}.

Indeed, the hardness of deciding whether objects of the kind of \eqref{eq:Mtn} are psd has been studied. Specifically:

 \begin{problem}
\label{pro} Given $ \{A_{j,j'} \in \Her_s \}_{j,j'=1}^{r}$  with $s,r\geq 7$, let 
$$M(n) :=  \sum_{j_1, \ldots , j_{n}=1}^{r}
A_{j_1,j_2} \otimes A_{j_2,j_3} \otimes \cdots \otimes 
A_{j_{n},j_1}.
$$
 Decide whether $M(n) \geq 0$ for all $n$.
\end{problem}

\begin{proposition}[\cite{De15}]
Problem \ref{pro} is undecidable.
\end{proposition}

This holds true even if all matrices $A_{j,j'}$ are diagonal and their entries rational.  
Variants of this problem are also undecidable   \cite{Kl14}, and deciding whether it is psd for a finite number of system sizes and with  open boundary conditions  is NP-complete  \cite{Kl14}.
See  also \cite{We14} for further perspectives on this problem.

\subsection{Example application: free spectrahedra}
\label{ssec:appl2}
Our second example application is in the area of convex optimization, where we find the same algebraic structures as the ones we have considered so far, albeit often studied from a  different angle.
 Namely, given a tuple  $(B_1,\ldots, B_{r})$ of Hermitian matrices, its associated 
\emph{spectrahedron} \cite{Bl13} is defined as 
$$
\mathrm{S}(B_1,\ldots, B_{r}) = \left\{(y_1, \ldots, y_{r}) \in \mathbb{R}^{r} \: \Bigg| \: \sum_{i=1}^{r} y_i B_i \geqslant 0\right\}, 
$$
where $\geqslant$ denotes positive semidefiniteness. 
Note that it is the intersection of an affine space with the set of psd matrices. 
Spectrahedra are precisely the feasible sets of semidefinite programimg (SDP). Characterizing which convex sets are spectrahedra is a main objective in the area of algebraic convexity \cite{Bl13}, and has direct implications for the many applications of semidefinite programming. 

Recently, a non-commutative generalization of spectrahedra has been proposed, called \emph{free spectrahedra} \cite{He13b}, defined as
\be
\mathrm{FS} (B_1,\ldots, B_{r}) = \bigcup_{s=1}^{\infty} \left\{  (A_1,  \ldots, A_d) \in \Her_s^{r} \: \Bigg| \: \sum_{i=1}^{r} A_i \otimes B_i \geqslant 0\right\}.
\label{eq:fs}
\ee
Thus, asking whether $\sum_{i=1}^{r} A_i \otimes B_i$ is psd is equivalent to asking whether $(A_1,\ldots, A_{r})$ is in the free spectrahedron defined by $(B_1,\ldots, B_{r})$. This is again a problem of the form~\eqref{eq:Cl}, with $n=2$. Surprisingly, many things about standard spectrahedra can be learned by examining their free versions. For example, the inclusion test of spectrahedra proposed in \cite{Be02} was   fully understood in \cite{He13} as an inclusion test of free spectrahedra, opening the way to analyzing the exactness of the method (see \cite{Fr16} and the references therein). Also, two minimal matrix tuples defining the same spectrahedron are unitarily equivalent if and only if they define the same free spectrahedron \cite{He13}. So the different free spectrahedra over a standard spectrahedron characterize equivalence classes of its defining matrix tuples. Applying these results often means checking whether a free spectrahedron contains a certain matrix tuple, i.e.\ whether $\sum_i A_i\otimes B_i$ is psd. Since the matrices might be of very large size, this is again a context in which our methods can be applied.

\subsection{Related work}

The methods  we use to solve the problems described above are fairly standard: a combination of techniques used for moment problems combined with results on sums of squares representations of positive polynomials. One might therefore expect there to exist a substantial amount of literature on the problem which we solve in this work, but this does not seem to be the case.

There is only one work that we are aware of. Lasserre~\cite{La11b} has investigated the smallest interval $[a_m,b_m]$ on which the support of a measure on $\R$ can be shown to be contained, given only its first $m$ moments. One can think of this as providing a solution to the problems discussed above in the case $p = \infty$. Lasserre found that the lower bound $a_m$ and upper bound $b_m$ are the optimal solutions of two simple semidefinite programs involving the moments.

\subsection{Overview}

The rest of this paper is structured as follows. 
In Section \ref{sec:pre}, we  characterize the distance of a matrix to the psd cone. 
In Section \ref{sec:idea}, we  provide upper and lower bounds to the distance to the psd cone by using a few moments of the matrix. 
In Section \ref{sec:algo}, we provide three methods to compute these bounds: 
the sos polynomial method, the Handelman method and the Chebyshev method. 
In Section \ref{sec:num} we analyse the numerical performance of these methods.
Finally, in Appendix \ref{vNa}, we will sketch the extension to von Neumann algebras.

\section{The distance to the psd cone}
\label{sec:pre}

We will start with some preliminaries (Section \ref{ssec:pre}) and then define the negative part function and the distance to the psd cone (Section \ref{ssec:np}).

\subsection{Preliminaries}
\label{ssec:pre}

Let us first fix some notation. 
For a matrix $M$, we denote its Hermitian conjugate by $M^*$.
$\Her_s$ denotes the set of Hermitian matrices of size $s\times s$. For $M\in\Her_s$, 
$M\geqslant 0$ denotes that $M$ is psd, and $M\leqslant 0$ that it is negative semidefinite (i.e.\ $-M\geqslant 0$).

We now state some basic facts about matrices and their norms. 
Consider $M\in \Her_s$ and  its spectral decomposition $M = U^*DU$, where $U$ is a unitary matrix, 
$D=\diag(\lambda_1,\ldots, \lambda_s)$, and $\{\lambda_1,\ldots, \lambda_s\}=:\spec (M)$ is the spectrum of $M$ (considered as a multiset). Any real-valued function $f$ defined on $\spec(M)$ can be defined on $M$ by setting $$f(M):=U^*f(D)\, U,$$ where $f(D)=\diag(f(\lambda_1),\ldots, f(\lambda_s))$.
For example, the absolute value $\vert M\vert$ of $M$ is defined this way.  
For $p\in [1,\infty)$, we define the Schatten $p$-norm of $M$ as
$$
\|M\|_p := \left(\ntr ( | M|^p ) \right)^{1/p},
$$
where taking the normalized trace instead of the usual trace introduces an additional factor of $s^{-1/p}$ with respect to the usual definition. This definition guarantees that $\|I\|_p = 1$, where $I$ is the identity matrix (of any size).
The case $p=2$ corresponds to the normalized  Frobenius norm, and the case $p=1$ to the normalized trace norm. 
Note also that if $p$ is even, the norm is easier to compute, since in this case the absolute value is superfluous, so that $\|M\|_p =  \left(\ntr (M^p) \right)^{1/p}$, which is simply the $p$-th root of the $p$-th moment of $M$. 
The operator norm of $M$ induced by the standard Euclidean norm on $\C^s$ is defined as
$$
\|M\|_\infty := \max_{\lambda_i\in \spec (M)} |\lambda_i| . 
$$
We have that $\|M\|_\infty = \lim_{p\to\infty} \|M\|_p$.

\begin{remark}\label{rem:scale} In the following we will often assume that $\Vert M\Vert_\infty\leq 1$, i.e.\ that the spectrum of $M$ lies in $[-1,1]$. This can clearly be achieved by a suitable scaling of $M$. Of course, since our main problem is that $\spec(M)$ cannot be computed, we cannot simply  scale by $\Vert M\Vert_\infty$.  But for $1\leq p\leq q\leq\infty$, we have
\[
	\Vert M\Vert_p\leq\Vert M\Vert_q \leq s^{\frac{1}{p}-\frac{1}{q}}\Vert M\Vert_p.
\]
So we can divide $M$ by $s^{1/p}\Vert M\Vert_p$ (for any $p\in\N$) in order to achieve $\|M\|_\infty \leq 1$. Moreover, since $s^{\frac{1}{\log s}} = e$, the norm $\Vert M\Vert_{\log s}$ is only a constant factor away from $\Vert M\Vert_\infty$.\footnote{This has been pointed out to us by Richard Kueng.}
\end{remark}

\subsection{The negative part function and the distance to the psd cone}
\label{ssec:np}

Given $p\in\N$, let us define the \emph{negative part function} $f_{p}$ as
\begin{equation}
\label{fp}
f_{p}(x) := \begin{cases}
0 & x\geq 0\\
\vert x\vert^p & x < 0. 
\end{cases}
\end{equation}
So, for example, $f_1(x)$ is the absolute value of the negative part of $x$.
For $M\in\Her_s$, we set 
$$
M_-:=f_{1}(M)\ \:\mbox{ and }\: M_+:=M+M_-,
$$
so that we obtain the natural definition $M=M_+-M_-$, where both the positive part $M_+$ and the negative part $M_-$ are psd. 
(To define functions of matrices we follow the standard procedure \cite{Hi08}. Namely, if the matrix is diagonalisable, as in our case $M = U D U^*$, then $f(M) = U f(D) U^*$, where $f(D) = \diag(f(\lambda_1),f(\lambda_2),\ldots )$, where $\{\lambda_i\}$  are the eigenvalues.)

Given a Hermitian matrix $M$ of size $s$, we are interested in its distance to the cone of psd matrices of size $s$, $\textrm{PSD}_s$, with respect to the Schatten $p$-norm, namely 
\be
d_p(M) := \inf_{N\in \textrm{PSD}_s} || M-N||_p .
\ee 
We now show that this is given by the negative part of $M$. That is,  the best psd approximation to a Hermitian matrix is given by the positive part of this matrix.

\begin{lemma}
\label{prop:psddist}
For $p\in[1,\infty)$ and $M\in\Her_s$, the matrix $M_+$ is the point in the psd cone that is closest to $M$, with respect to any Schatten $p$-norm. That is, 
\[
	d_p(M) = \Vert M_-\Vert_p = \ntr(f_{p}(M))^{1/p} .
\]
\end{lemma}

\begin{proof} Clearly $M_+\geqslant 0$, so that the distance can at most be $\Vert M_-\Vert_p$. Now, for any matrix $A\in\Her_s$, denote by $\sigma(A)$ the diagonal matrix with the eigenvalues of $A$ on the diagonal, in decreasing order. It follows from \cite{Bh91} (IV.62) that
\begin{equation}
\label{eq:major}
	\Vert P-M\Vert_p \geq \Vert \sigma(P)-\sigma(M)\Vert_p
\end{equation}
holds for all $P\in\Her_s$ and all $p$. If $P\geqslant 0$, then clearly $$\Vert P-M\Vert_p \geq \Vert \sigma(P)-\sigma(M)\Vert_p \geq \left(\frac{1}{s}\sum_{\lambda\in\spec(M), \lambda<0} \vert \lambda\vert^p\right)^{1/p}=\Vert M_-\Vert_p.$$ This proves the claim.\qed
\end{proof}

\begin{remark}
Note that $\Vert \sigma(P)-\sigma(M)\Vert_p$ is precisely the $p$-Wasserstein distance \cite{Vi09} between the spectral measures  $s^{-1} \sum_i \delta_{\sigma(P)_{ii}}$ and $s^{-1} \sum_i \delta_{\sigma(M)_{ii}}$ of $P$ and $M$. 
So $d_p(M)$ is in fact also the $p$-Wasserstein distance of the spectral measure of $M$ to the cone  of probability measures supported on $[0,\infty)$.
\demo
\end{remark}

\section{Bounds on the distance}\label{sec:idea}

We will now bound  the negative part function $f_p$ by polynomials (Section \ref{ssec:fppoly}) 
and then show the asymptotic optimality of our bounds (Section \ref{ssec:opt}).

\subsection{Bounding the negative part function by polynomials}
\label{ssec:fppoly}

 Clearly,  $M$ is psd if and only if its distance to the psd cone is zero, $d_p(M) = 0$, and $d_p(M)$ is a measure of how far $M$ differs from being psd. Thus, ideally, we would like to compute $d_p(M)$. Since we only have access to a few moments of $M$, we will use them to estimate $d_p(M)$ as accurately as possible. 
 
We start by showing that
$d_p(M)$ can be approximated arbitrarily well by the trace of polynomial expressions in $M$. 
For $q\in\R[x]$, we write $f_{p}\leqslant q$ if this holds pointwise on the interval $[-1,1].$

\begin{lemma}\label{lem:app}
Suppose  $M\in\Her_s$ with $\Vert M\Vert_\infty\leq 1$ and $p\in\N$. Then 
\[
	d_p(M)^p = \inf_{f_{p}\leqslant q\in\R[x]}\ntr ( q(M) )= \sup_{f_{p}\geqslant q\in\R[x]}\ntr ( q(M) ).
\]
\end{lemma}
\begin{proof}
From $f_{p}\leqslant q$ we obtain that $f_{p}(M)\leqslant q(M)$, thus $\ntr(f_{p}(M))\leq \ntr(q(M))$ and finally  $$d_p(M)^p=\ntr(f_{p}(M))\leq \ntr(q(M)).$$ Conversely, by standard Weierstrass approximation for continuous functions on compact sets, for any $\epsilon >0$ there exists a polynomial $q\in\R[x]$ with $$f_{p}\leqslant q\leqslant f_{p}+\epsilon.$$ The same  argument as above then shows $$\ntr(q(M))\leq d_p(M)^p + \epsilon.$$ This proves the first equation, and the second follows in the same fashion. \qed
\end{proof}

\begin{figure}[htb]
\includegraphics[scale=0.6]{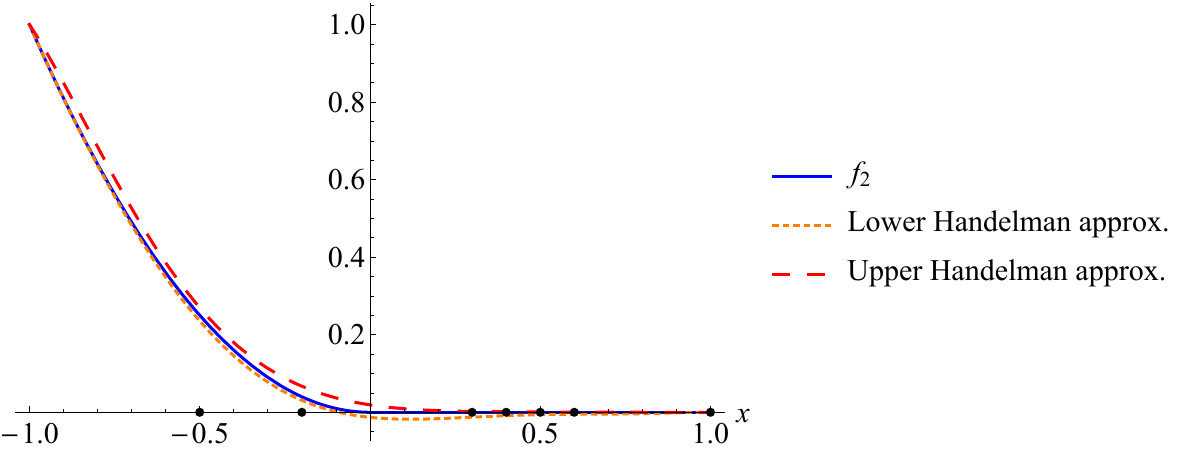}
\caption{Polynomials of degree 7 that approximate $f_2$ from below and from above, together with the spectrum of $M$ (black dots). The polynomials have been obtained with the Handelman method, to be described in  Section \ref{sec:method2}. }\label{fig:app1}
\end{figure}

Thus,  any polynomials $q_1$ and $q_2$ such that 
  $q_1\leqslant f_p\leqslant q_2$ give the bounds
\[
	\ntr(q_1(M))\leq d_p(M)^p\leq \ntr(q_2(M)).
\]
In particular, this means that if  $\ntr(q_1(M))>0$  then $M$ is not psd. 
The quality of the bounds depends on the quality of the approximation of $f_p$ by $q_1,q_2$ on $\spec(M)$, or more generally on $[-1,1]$ --- see Figure \ref{fig:app1} for an example.

\begin{figure}[htb]
\includegraphics[scale=0.6]{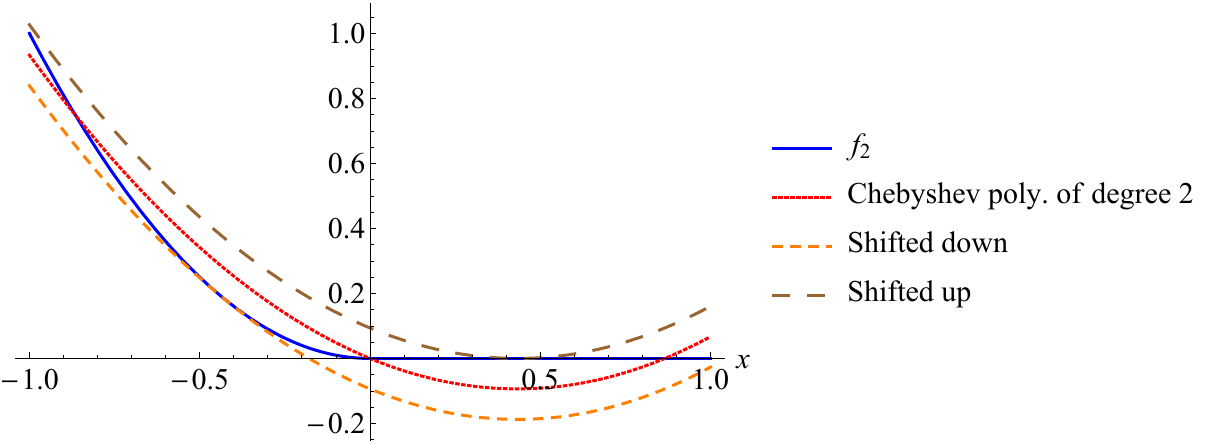}
\caption{A polynomial $q$ which approximates $f_2$ but does not satisfy $q\leqslant f_2$ or $q\geqslant f_2$. 
The polynomial is a Chebyshev polynomial of degree 2, obtained with the Chebyshev method (Section \ref{sec:method1}). 
Then it is shifted up and down by $|| q-f_2||_\infty$. 
 }
 \label{fig:app2}
\end{figure}

More generally, given a polynomial $q$ which does not satisfy $q\leqslant f_p$ or $q\geqslant f_p$, what we can say is that 
\be
	q - \Vert (q - f_p)_+ \Vert_\infty \leqslant f_p \leqslant q + \Vert (q - f_p)_- \Vert_\infty ,
	\label{eq:shift}
\ee
where we use notation for positive and negative part as in the matrix case, and similarly $\Vert g\Vert_\infty :=\sup_{x\in[-1,1]} | g(x)|$. Here, the function on the left hand side corresponds to shifting $q$ by an additive constant until it is below $f_p$, and similarly for the right hand side --- see Figure \ref{fig:app2} for an example. 
This leads to the following result. 

\begin{theorem}\label{thm:main}
Let $M\in\Her_s$ with $\Vert M\Vert_\infty\leq 1$, $p\in\N$, as well as $q\in\R[x]$. Then
\be
\label{eq:bound}
	\qquad \ntr(q(M)) - \Vert (q-f_p)_+ \Vert_\infty \:\leq\: d_p(M)^p \:\leq\: \ntr(q(M)) +\Vert (q-f_p)_- \Vert_\infty.
\ee
In particular:
\begin{itemize}
\item[(i)] If $f_p\leqslant q$, then $d_p(M)\leq \ntr(q(M))^{1/p}$. 
\item[(ii)] If $\ntr(q(M)) > \Vert (q-f_p)_+ \Vert_\infty$, then $M$ is not psd.
\item[(iii)] If $q\leqslant f_p$ and $\ntr(q(M))>0$, then $M$ is not psd.
\end{itemize}
\end{theorem}

\begin{remark}
We can also try to prove that a matrix is psd with the same approach, although this seems to be much harder in practice. Note the following: if for $M\in\Her_s$ we find $d_p(M)\leq \epsilon$, then $d_\infty(M)\leq s^{1/p} \eps$, and therefore $M+s^{1/p}\epsilon I_s\geqslant 0$. So if we first replace $M$ by
\[
	M_\epsilon := M-s^{1/p}\epsilon I_s,
\]
and then  show that $d_p(M_\epsilon)\leq\epsilon$, then we have proven that $M\geqslant 0$. By Theorem \ref{thm:main}, this can be achieved by finding $q\in\R[x]$ with
\be
\label{eq:pert}
\ntr(q(M_\epsilon)) +\Vert (q-f_p)_- \Vert_\infty\leq \epsilon^p.
\ee
Note that the $k$-th moment of $M_\epsilon$ can be computed from the moments of $M$ of order $\leq k$. Further, if $M$ is strictly positive definite, then this strategy does indeed work: there is $\epsilon>0$ with $M_\epsilon \geqslant 0$, i.e. $d_p(M_\epsilon)=0$. In view of Lemma \ref{lem:app}, we can also find some $q$ to make (\ref{eq:pert}) hold, so that this method can indeed detect that $M$ is psd. However, we might need to test for both very small $\epsilon$ and polynomials $q$ of very high degree before obtaining a positive answer.
\demo
\end{remark}

\subsection{Optimality of the bounds}
\label{ssec:opt}

We now show that given the first $m$ moments of $M$, 
the best bounds to  the distance to the psd cone which are independent of the size of $M$ are given by a polynomial approximation. 
More precisely, given $M\in\Her_s$ and $p,m\in\N$, define 
\begin{align}
\begin{split}
\label{eq:pn}
d_{p,m}^+(M) &:= \inf_{f_{p}\leqslant q,\: \deg(q)\leq m}\ntr( q(M) )^{1/p} \\ d_{p,m}^-(M) &:= \sup_{q\leqslant f_{p},\: \deg(q)\leq m}\ntr( q(M) )^{1/p} .
\end{split}
\end{align}
Clearly,  these numbers depend only on the first $m$ moments of $M$, and 
they lower and upper bound $d_p(M)$, 
\be
\label{eq:opt}
0\leq d_{p,m}^-(M)\leq d_p(M)\leq d_{p,m}^+(M). 
\ee
We now show that these are the optimal upper and lower bounds to $d_p(M)$ which can be obtained from the first $m$ normalized moments of $M$ and which are independent of the size of $M$. We thus call them \emph{asymptotically optimal}.

The following result is a variation on a classical result from the theory of moment problems~\cite[Theorem~4.1(a)]{Kr77}, for which we offer a self-contained proof.

\begin{theorem}
For any matrix $M\in\Her_s$ with $\Vert M\Vert_\infty\leq 1$ and any $m\in\N$, 
\begin{itemize}
\item For every $\epsilon > 0$, there are $N_1,N_2\in\Her_{t}$ (for suitably large $t$) such that 
 $$
 \vert \ntr(M^k) - \ntr(N_i^k)\vert \leq\epsilon \quad \textrm{for } k=1,\ldots, m ,
 $$
 and for which 
$$
d_p(N_1)\geq d_{p,m}^+(M), \quad 
d_p(N_2) \leq d_{p,m}^-(M).
$$
\item 
There are operators $N_1$ and $N_2$ in a finite von Neumann algebra $\mathcal{N}$ such that  
$$
\ntr(M^k) = \ntr(N_i^k) \quad \textrm{for } k=1,\ldots, m ,
$$
and which saturate the bounds, 
$$
d_p(N_1)=d_{p,m}^+(M)
, \quad 
d_p(N_2) = d_{p,m}^-(M).
$$
\end{itemize}
\label{thm:optimality}
\end{theorem}

Note that in the first case the moments are approximately reproduced but $N_i$ has finite size --- we will see an example thereof in Example \ref{ex}.
In the second case, the moments are exactly reproduced but the size of $N_i$ may need to be  infinite.

\begin{proof} We only construct $N_1$, since  $N_2$ is obtained similarly. 

Consider the  linear functional 
\begin{align*}
\varphi\colon \R[x]_{\leq m}&\to \R \\ q&\mapsto \ntr(q(M))
\end{align*} 
 on the space of polynomials of degree at most $m$, which clearly maps polynomials nonnegative on $[-1,1]$ to nonnegative numbers.  Let us define the real vector space
\[
	V:=\{\,q+rf_p \mid q\in \mathbb{R}[x]_{\leq m}, r\in \mathbb{R} \, \} .
\]
We extend $\varphi$ to a linear functional $\psi$ on $V$ by setting  $$\psi(f_p):=\inf_{f_p\leqslant q\in\R[x]_{\leq m}} \varphi(q).$$
 We  claim that $\psi$ still maps nonnegative functions in $V$ to nonnegative numbers. So let $h+rf_p\geqslant 0$ for some $h\in\R[x]_{\leq m}$ and $r\in\R$. The case $r = 0$ is clear, so assume $r<0$. Then $-\frac{1}{r}h\geqslant f_p$, and thus
\begin{align*}
	\psi(h+rf_p) & = \varphi(h) + r \inf_{f_p\leqslant q} \varphi(q) \\ 
	& \geq \varphi(h) + r \varphi\left(-\frac{1}{r}h\right) \\ 
	& =\varphi(h)-\varphi(h)=0.
\end{align*}
If $r>0$ instead, then $f_p\leqslant q$ implies $0\leqslant h+rf_p\leqslant h+rq$, and thus
\[
	\varphi(h)+r\varphi(q)=\varphi(h+rq)\geq 0.
\]
Passing to the infimum over these $q$ proves the statement.

Since $\R[x]_{\leq m}$ already contains an interior point of the convex cone of nonnegative continuous functions on $[-1,1]$ (such as the constant function $1$), we can further extend $\psi$ to a positive linear functional $\Psi$ on the whole of $\mathcal C([-1,1]),$ using the Riesz Extension Theorem \cite{Ri23}. By the Riesz Representation Theorem \cite{Ru91}, there exists a positive measure $\mu$ on $[-1,1]$ such that 
$$ \Psi(f)=\int f\, d\mu$$ 
for all $f\in\mathcal{C}([-1,1])$. From $\mu([-1,1])=\int 1\, d\mu=\Psi(1)=\varphi(1)=1$,   we see that $\mu$ is a probability measure.
We now take $\mathcal{N} := L^\infty([-1,+1],\mu)$, equipped with integration against $\mu$ as a finite normalized trace, and define $N_1$ to be the multiplication operator by the identity function, $N_1 : f \mapsto xf$.

So for $k=0,\ldots, m$, the $k$-th moment of $N_1$ is given by
\[
	\int x^kd\mu=\Psi(x^k)=\varphi(x^k)=\ntr(M^k),
\]
and we also have
\[
	d_p(N_1) = \int f_p\, d\mu = \Psi(f_p) = \psi(f_p) = \inf_{f_p\leqslant q \in \R[x]_{\leq m}}\varphi(q) = d^+_{p,m}(M)^p,
\]
which establishes the first claim.

Concerning the realization by finite-dimensional matrices, we use the well-known fact that each probability measure on $[-1,1]$  can be approximated arbitrarily well by uniform atomic measures with respect to the weak-$^*$ topology. Concretely, we can approximate $\mu$ by the uniform atomic measure $\nu_t := \frac{1}{t}\cdot\sum_{i=1}^t \delta_{a_i}$, where the $a_i$ are the right $t$-quantiles,
\[
	a_i := \inf \left\{\: r \in [-1,1] \:\bigg|\: \mu([r,1]) \geq \frac{i}{t} \:\right\}.
\]
This choice guarantees that the cumulative distribution function of $\nu_t$ dominates the one of $\mu_t$. Therefore the expectation value of $\mu$ is not smaller than that of $\nu_t$  on any monotonically nonincreasing function. In particular, we have
\[
	\int f_p \, d\nu_t \geq \int f_p \, d\mu .
\]
Furthermore, since $\nu_t$ weak-$^*$ converges\footnote{One way to see this is by the Portmanteau theorem: the cumulative distribution functions of $\mu$ and $\nu_t$ differ by at most $t^{-1}$ at every point, and therefore we have (even uniform) convergence as $t\to\infty$, which implies weak-$^*$ convergence $\nu_t \to \mu$.} to $\mu$ as $t\to\infty$, we can moreover choose $t$ large enough such that
\[
	\left\vert \int x^k d\mu -\int x^k d\nu_t \right\vert \leq  \epsilon\quad  \mbox{ for } k=0,\ldots, m.
\]
For $N_1 := \diag(a_1,\ldots, a_t)\in\Her_t$,
we then have $\Vert N_1\Vert_\infty\leq 1$ and 
\[
	\ntr(N_1^k)= \frac{1}{t}\sum_{i=1}^t a_i^k=\int x^kd\nu_t
\]
as well as
\[
	\int f_p\, d\nu_t=\frac{1}{t}\sum_{a_i<0}\vert a_i\vert^p = \Vert {N_{1}}_-\Vert_p^p= d_p(N_1)^p,
\]
which gives the desired bounds. This altogether finishes the proof. \qed
\end{proof}

\begin{example}\label{ex}
We show that $\epsilon = 0$ can in general not be achieved in the second part of Theorem \ref{thm:optimality}.
Taking $m = 2$, let us consider the matrix\footnote{Although $M$ is only of size $2\times 2$, one can clearly achieve the same moments on larger matrices by simply repeating the eigenvalues.}
\[
	M = \diag\left[c + \sqrt{c(1-c)},\: c - \sqrt{c(1-c)}\right],
\]
for $c \in (0,1/2)$. Since the lower right entry is negative, $M$ is not psd. Its first moment is $c$, the second moment is $\frac{1}{2}\left(2c^2 + 2c(1-c)\right) = c$, equal to the first moment. Looking for a probability measure $\mu$ on $[0,1]$ with these moments, we must have $\E_\mu[x(1-x)] = 0$, which implies that $\mu$ must be supported on $\{0,1\}$ only. Since $\mu = (1-c) \delta_0 + c \delta_1$ does indeed have these moments, we conclude that it is the unique measure on $[0,1]$ with these moments; and the fact that such a $\mu$ exists implies $d_{p,2}^-(M) = 0$, irrespectively of the value of $p$. However, as soon as $c$ is irrational, the measure $(1-c)\delta_0 + c\delta_1$ is not of the form $s^{-1}\sum_{i=1}^s \delta_{\lambda_i}$ for any finite sequence $(\lambda_1,\ldots,\lambda_s)$. In particular, there is no psd matrix of finite size with the same moments as $M$, and $\epsilon = 0$ cannot be achieved in Theorem \ref{thm:optimality}. By a standard compactness argument, this also implies that one must take $t \to \infty$ as $\epsilon \to 0$ in the theorem.
\demo
\end{example}

\begin{remark}
An interesting question\footnote{communicated to us by Boaz Barak.} is how close our bounds $d_{p,m}^+(M)$ and $d_{p,m}^-(M)$ are guaranteed to be to the actual value $d_{p,m}(M)$. In the worst case, the actual value will coincide with one of the two bounds, in which case the other bound differs from the actual value by
\[
d_{p,m}^+(M) - d_{p,m}^-(M) 
\leq  12\frac{p}{m},
\]
which follows from \eqref{eq:bound2} and \eqref{cheb_bound}. 
\end{remark}

\section{Algorithms}
\label{sec:algo}

We now present our numerical methods to compute lower and upper bounds to the distance to the psd cone, in order of decreasing accuracy and complexity: 
the sos polynomial method (Section \ref{sec:method3}),
the Handelman method (Section \ref{sec:method2})
and the Chebyshev method (Section \ref{sec:method1}). 
The sos polynomial method involves solving a semidefinite program,
the Handelman method involves solving a linear program, 
and the Chebyshev method does not require any optimization.
We will compare the numerical performance of three methods in  Section \ref{sec:num}.

Throughout, we fix nonzero 
 $p\in\N$.

\subsection{The sos polynomial method}
\label{sec:method3} 

The sos polynomial method solves the optimization problems of Eq.\ \eqref{eq:pn} exactly\footnote{This is not to be confused with the sos polynomial method of \cite{De13c}, which is a semidefinite program that computes $\min_p \Vert M - p(M) \Vert_1$, where $p$ is a sos polynomial of given degree $m$. The goal of the method of \cite{De13c} is to approximate $M$ as well as possible with a sos polynomial (as this provides a purification), which is possible only if $M$ is psd. Note moreover that  $\Vert M - p(M) \Vert_1$ cannot be computed from the moments of $M$.
},
and thereby computes $d_{p,m}^+(M)$ and $d_{p,m}^-(M)$. 
We start by explaining how to compute the upper bound $d_{p,m}^+(M)$ via a semidefinite program.

To be able to talk about polynomials only, we first split the condition $f_p\leqslant q$  into two parts:
\begin{equation}
\label{eq:condq}
	0\leqslant q(-x)-x^p\ \mbox{ and }\  0\leqslant q(x),\ \mbox{ both for all } x\in [0,1].
\end{equation}
 Now note that any polynomial of the form
\begin{equation}
\label{eq:posansatz}
	\sigma_0+\sigma_1x+\sigma_2(1-x)+\sigma_3x(1-x),
\end{equation}
where the $\sigma_i$ are sums of squares of polynomials, is  nonnegative on $[0,1]$. In fact, the converse is true as well:

\begin{theorem}[\cite{Ku02,Ku05}]\label{thm:pos}
If $q\in\R[x]_{\leq m}$ is nonnegative on $[0,1]$, then there exist sums of squares $\sigma_0,\sigma_1,\sigma_2,\sigma_3\in\R[x]$ such that
\[
	q(x)=\sigma_0+\sigma_1x+\sigma_2(1-x)+\sigma_3x(1-x)
\]
where the degree of each $\sigma_i$ can be chosen such that each summand has degree at most $\leq m$.
\end{theorem}

So assume that we find such representations for both polynomials in (\ref{eq:condq}),

\begin{align*}
	q(-x)-x^p & = \sigma_0+\sigma_1x+\sigma_2(1-x)+\sigma_3x(1-x), \\
	q(x) & = \tau_0+\tau_1x+\tau_2(1-x)+\tau_3x(1-x), 
\end{align*}
with sums of squares (sos) $\sigma_i,\tau_i\in\R[x]$.  Then we have clearly ensured $f_p\leqslant q$. This can be rewritten as 
\begin{equation}
\label{eq:sdp}
\begin{split}
	q(x) &= (-x)^p  + \tilde\sigma_0 - \tilde\sigma_1 x + \tilde\sigma_2 (1+x) -\tilde\sigma_3 x (1+x) \\
	& = \tau_0+ \tau_1 x + \tau_2 (1-x) +\tau_3 x (1-x),
\end{split}
\end{equation}
where the $\tilde\sigma_i (x) := \sigma_i(-x)$ are again sos.

Now assume that every term in (\ref{eq:sdp}) has degree $\leq m$. This imposes obvious degree bounds on the sums of squares $\tilde\sigma_i,\tau_i$, namely $\deg(\tilde\sigma_0) \leq m$, and $\deg(\tilde\sigma_1), \deg(\tilde\sigma_2) \leq m-1$ as well as $\deg(\tilde\sigma_3) \leq m-2$, and analogously for the $\tau_i$, and note that every sos polynomial must have an even degree. It is easy to see that every sum of squares can be written as
\be
\begin{split}
&\tilde\sigma_i = (1,x,\ldots,x^{l_i}) \, S_i\, (1,x,\ldots,x^{l_i})^t, \quad S_i\geq 0\\
&\tau_i = (1,x,\ldots,x^{l_i}) \, T_i\, (1,x,\ldots,x^{l_i})^t, \quad T_i\geq 0
\end{split}
\label{eq:sdp2}
\ee
where $l_i = \textrm{deg}(\tilde\sigma_i)/2$ and similarly for $\tau_i$. 
 Writing each $\sigma_i,\tilde\tau_i$ in such a way, using matrices with unknown entries, and then comparing coefficients with respect to $x$ in (\ref{eq:sdp}), leads to the problem of finding psd matrices with certain linear constraints on the entries. Any solution to this problem will provide a polynomial $q\in\R[x]_{\leq m}$ with $f_p\leqslant q.$ Among all of them, we want to minimize $\ntr(q(M))$, which is a linear function in the entries of the unknown matrices, having the moments of $M$ as coefficients.  Optimizing a linear function over an affine section of a cone of psd matrices is  known as semidefinite programming, for which there exist efficient algorithms.

The derivation of the lower bound $d_{p,m}^-(M)$ is entirely parallel, except for the fact that  in \eqref{eq:condq} the two inequalities are reversed.  This implies that
\begin{align*}
	-q(-x)+x^p & = \sigma_0+\sigma_1x+\sigma_2(1-x)+\sigma_3x(1-x), \\
	-q(x) & = \tau_0+\tau_1x+\tau_2(1-x)+\tau_3x(1-x), 
\end{align*}
where $\sigma_i,\tau_i$ are sos, and thus 
\begin{equation}
\label{eq:sdplb}
\begin{split}
	q(x)= & ( -x)^p - \tilde\sigma_0 + \tilde\sigma_1 x - \tilde\sigma_2 (1+x) +\tilde\sigma_3 x (1+x) \\
	= & - \tau_0- \tau_1 x - \tau_2 (1-x) -\tau_3 x (1-x).
\end{split}
\end{equation}

In summary, we have obtained:

\begin{proposition}[Sos polynomial method]
The sos polynomial method at level $m$ computes the upper bound $d_{p,m}^+(M)^p$ (defined in \eqref{eq:pn}) by solving the semidefinite program
\begin{equation}
\begin{split}
&\min \: \ntr(q(M)) \\
&\mathrm{subject\: to\: \: Eq.\:} \eqref{eq:sdp}\\
&\quad\qquad \mathrm{and\: \: Eq.\:}  \eqref{eq:sdp2}
\end{split}
\end{equation}
It also computes the lower bound $d_{p,m}^-(M)^p$ (defined in \eqref{eq:pn}) by solving the semidefinite program
\begin{equation}
\begin{split}
&\max \: \ntr(q(M)) \\
&\mathrm{subject\: to\: \: Eq.\:} \eqref{eq:sdplb}\\
&\quad\qquad \mathrm{and\: \: Eq.\:}  \eqref{eq:sdp2}
\end{split}
\end{equation}
\end{proposition}

\begin{figure}[htb]
\includegraphics[scale=0.6]{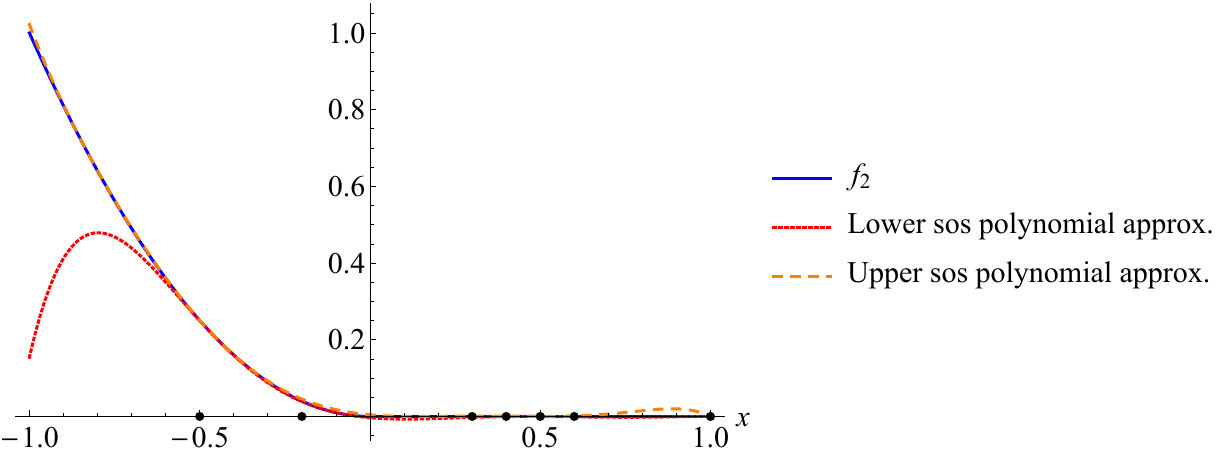}
\caption{
Lower and upper sos polynomial approximation of degree 7 for $p=2$ and for the matrix  $M$ whose spectrum is shown in black dots. 
} 
\label{fig:app4}
\end{figure}

As an example, Figure \ref{fig:app4} shows the sos polynomial approximation of degree $7$ obtained for a matrix $M$ with the indicated spectrum. We will discuss numerical results more systematically in  Section \ref{sec:num}. 

\begin{remark}
For computational purposes, it may be better to use the ansatz
\be
\label{eq:simp}
	\sigma_0 +\sigma_1x(1-x)
\ee
for a polynomial nonnegative on $[0,1]$, instead of \eqref{eq:posansatz}. The advantage is that this reduces the number of sums of squares from $4$ to $2$. An analogue of  Theorem \ref{thm:pos} still holds \cite{Ku05},  but with slightly weaker degree bounds: a degree $m$ polynomial nonnegative on $[0,1]$ has a representation $\sigma_0 +\sigma_1x(1-x)$ with $\deg(\sigma_0)\leq m+1$ and $\deg(\sigma_1)\leq m-1$ only.\footnote{This can also be deduced directly from Theorem \ref{thm:pos}, using $x=x^2+x(1-x)$ and $1-x=(1-x)^2+x(1-x)$.} So if we set up our optimization problem as above, but with the simpler representation (\ref{eq:simp}) where we demand $\deg(\sigma_0)\leq m$ and $\deg(\sigma_1)\leq m-2$ (since moments of $M$ are available only up to order $m$), then we will obtain a bound on $d_p(M)$ which lies between $d_{p,m}^+(M)^p$ and $d_{p,m-1}^+(M)^p$.
\demo
\end{remark}

\subsection{The Handelman method}
\label{sec:method2}
The  Handelman method relaxes the optimization of Eq.\ \eqref{eq:pn} by using another ansatz for nonnegative polynomials on $[0,1]$. This results in a linear optimization problem, which can usually be solved much faster than a semidefinite problem. 

We start by splitting the condition $f_p\leqslant q$ as in \eqref{eq:condq}. This time note that any polynomial of the form
\[
	\sum_{\alpha\in\N^2} b_\alpha x^{\alpha_1}(1-x)^{\alpha_2}
\]
with $b_\alpha\geq 0$ (and only finitely many different from 0) is nonnegative on $[0,1]$. So if we find coefficients $b_\alpha, c_\alpha\geq 0$ with 
\begin{align*}
	q(-x)-x^p & = \sum_{\alpha\in\N^2} b_\alpha x^{\alpha_1}(1-x)^{\alpha_2}, \\[4pt]
	q(x) & = \sum_{\alpha\in\N^2} c_\alpha x^{\alpha_1}(1-x)^{\alpha_2},
\end{align*}
then we can be sure to have $f_p\leqslant q$. 
Assume that the degree of the polynomials is a priori bounded by $m$. Then comparing coefficients with respect to $x$ yields a finite system of linear equations:
\be
\begin{split}
q(x) &= (-x)^p + \sum_{
|\alpha|\leq m
} b_\alpha (-x)^{\alpha_1}(1+x)^{\alpha_2}  \\
&= 
 \sum_{
 |\alpha|\leq m
 } c_\alpha x^{\alpha_1}(1-x)^{\alpha_2}, 
 \end{split}
 \label{eq:q}
\ee
where $|\alpha| = \alpha_1+\alpha_2$.
We look for solutions under the constraint that $(b_\alpha,c_\alpha)_{|\alpha|\leq m}$, and 
among all these solutions, we look for the one that minimizes the quantity $\ntr(q(M))$. 
This is precisely a linear optimization problem, where information about our matrix $M$ enters through the objective function. 

The derivation of the lower bound of $d_{p,m}^-(M)$ is analogous, except that in this case 
$- q(-x)+x^p $ and $-q(x)$ must be nonnegative polynomials on $[0,1]$. This leads to 
\be
\begin{split}
q(x)&= (-x)^p - \sum_{
|\alpha|\leq m
} b_\alpha (-x)^{\alpha_1}(1+x)^{\alpha_2}  \\
&= 
 \sum_{
 |\alpha|\leq m
 } - c_\alpha x^{\alpha_1}(1-x)^{\alpha_2}, 
\end{split}
 \label{eq:qlb}
\ee

In summary we have obtained: 

\begin{proposition}[Handelman method]
\label{met:handelman}
The Handelman method at level $m$ computes an upper bound of $d_{p,m}^+(M)$ by solving the linear program
\begin{equation}
\begin{split}
&\min \: \ntr(q(M)) \\
&\mathrm{subject\: to\: \: Eq.\:} \eqref{eq:q}\\
&\quad\qquad \mathrm{and\:} b_\alpha\geq 0, c_\alpha\geq 0
\end{split}
\end{equation}
It also computes a lower bound of $d_{p,m}^-(M)$ by solving the linear program
\begin{equation}
\begin{split}
&\max \: \ntr(q(M)) \\
&\mathrm{subject\: to\: \: Eq.\:} \eqref{eq:qlb}\\
&\quad\qquad \mathrm{and\:} b_\alpha\geq 0, c_\alpha\geq 0
\end{split}
\end{equation}
\end{proposition}

Note that linear optimization problems are easy to solve algorithmically (for example, interior point methods have polynomial complexity, but the simplex algorithm and its variants often work best in practice). 
See Figure \ref{fig:app1} above for upper and lower Handelman approximations of $f_2$ for a given matrix $M$.

Although this method computes only a  relaxation of $d_{p,m}^+(M)$ (and analogously for $d_{p,m}^-(M)$), the following special case of Handelman's Positivstellensatz for polytopes \cite{Ha88} ensures that these relaxations converge to the exact value $d_p(M)$ in the limit $m\to\infty$:

\begin{theorem}[Handelman]
If $q\in\R[x]$ is strictly positive on $[0,1]$, then $$q(x)=\sum_{\alpha\in\N^2} a_\alpha x^{\alpha_1}(1-x)^{\alpha_2}$$ for certain $a_\alpha\geq 0,$ only finitely many different from 0.
\end{theorem}

Note that this result leads directly to the standard solution of the Hausdorff moment problem in terms of complete monotonicity~\cite[Theorem~1.5]{Sh70}. We now have:

\begin{corollary}
Let $\tilde d_{p,m}^+(M)$ denote the $p$-th root of the optimal value of the linear program described in Method \ref{met:handelman}.  Then we have $d_p(M)\leq d_{p,m}^+(M)\leq \tilde d_{p,m}^+(M)$, and
\[
	\tilde d_{p,m}^+(M)\stackrel{m\to\infty}{\searrow} d_p(M).
\]
\end{corollary}

\begin{proof} Every feasible point in the program leads to a polynomial $q\in\R[x]_{\leq m}$ with $f_p\leqslant q$.  This proves  $d_{p,m}^+(M)\leq \tilde d_{p,m}^+(M)$. Now if $q\in \R[x]$ fulfills $f_p\leqslant q,$ then for any $\epsilon >0$ there is some $m\in\N,$ such that $q+\epsilon$ corresponds to a feasible point in the respective program, by Handelman's theorem. This proves the claim. \qed
\end{proof}


\subsection{The Chebyshev method}
\label{sec:method1}

The Chebyshev method  chooses $q_m$ as the Chebyshev polynomial of degree $m$ that best approximates $f_p$, and uses \eqref{eq:bound} to derive bounds on $d_p(M)$. This has the advantage of not involving any optimization at all, i.e.\ one need only compute the Chebyshev polynomials $q_m$ once, and the method can be applied to any matrix $M$.

Ideally, the best bounds of Eq.\ \eqref{eq:bound} are given by the polynomial $q$ of degree $m$
 that minimizes $||q_m-f_p||_\infty$, which we call $q^*_m$. 
By Jackson's theorem (see e.g.\ \cite[Theorem 1.4]{Ri81}) we have that 
\begin{equation}
\label{optconv}
	\Vert f-q^*_m\Vert_\infty \leq 6\, \omega(f,1/m).
\end{equation}
where $\omega(f,\delta)$ denotes the modulus of continuity of a uniformly continuous function $f$, defined as
\[
	\omega(f,\delta) = \sup_{\begin{subarray}{l}x_1,x_2\in[-1,1]\\
	|x_1-x_2|\leq\delta \end{subarray}}
	|f(x_1)-f(x_2)|.
\]
For the negative part function $f_p$, we thus obtain that 
\begin{equation}
\label{cheb_bound}
\|f_p - q^*_m\| _{\infty}\leq 6 \left(1-\left(1-m^{-1}\right)^p\right) \leq 6 \frac{p}{m},
\end{equation}
where the second estimate is by Bernoulli's inequality, which is tight up to an error of $O(m^{-2})$.
 
\begin{figure}[htb]\centering
\includegraphics[scale=0.6]{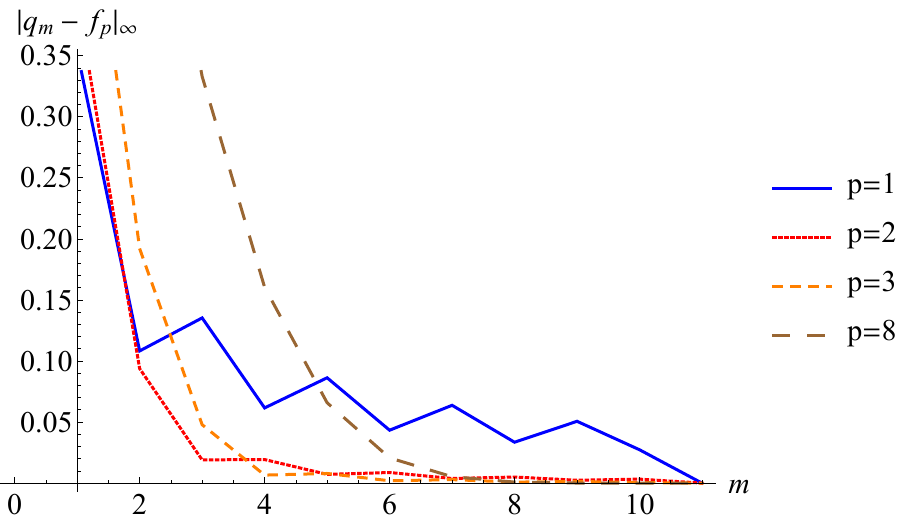}
\caption{$|| q_m-f_p||_\infty$ (obtained numerically) for the Chebyshev polynomial $q_m$ of \eqref{eq:qm} as a function of $m$ for several $p$. 
}\label{fig:InfNormCheby}
\end{figure}

However, it is in general not possible to find $q_m^*$ (called the minimax approximation)  analytically. Instead, Chebyshev polynomials provide a good proxy: they are close to the minimax approximation, and are  straightforward to compute explicitly \cite{Ri81}. To obtain an analytical upper bound of the additional error, it is known \cite{Gi07} that the  Chebyshev approximation $q_m$ of a continuous function $f$ satisfies that
$$
\Vert q_m-f\Vert_\infty\leq \, C \omega(f,1/m) \log m \,
$$
for some constant $C$. For our negative part function $f_p$, it thus holds that 
\be
|| q_m-f_p||_\infty \leq C m^{-1} \log m .
\label{eq:infche}
\ee
Numerically we can see the behaviour of  $||q_m-f_p||_\infty$ as a function of $m$ and $p$ in Figure \ref{fig:InfNormCheby}.

Let us now construct the Chebyshev interpolating polynomial of degree $m$, $q_m$, explicitly. 
It interpolates $f_p$ at the Chebyshev nodes  
$$
 a_k=\cos(\pi (k+1/2)/(m+1)) \quad \mbox{ for } k=0,1,\ldots,m, 
 $$
so that   $q_m$  is given by 
$$
q_m = \sum_{k=0}^m f_p(a_k) \ell_k ,
$$
 where the $\ell_k$ are the Lagrange interpolation polynomials 
$$
\ell_k (x)=\frac{1}{\prod_{i\neq k}(a_k-a_i)} \prod_{\begin{subarray}{c}
i\neq k
\end{subarray}}(x-a_i).
$$
One can also express $q_m$ in terms of the Chebyshev polynomials of the first kind, which are defined as 
$
t_k(x) = \cos(k\arccos(x))$ for $|x|\leq 1.$
In this basis, $q_m$ takes the form
\be
q_m = {\sum_{k=0}^m} c_k t_k -c_0/2, \qquad \textrm{with }\: c_k = \frac{2}{m+1}\sum_{j=0}^m 
f_p(a_j) t_k(a_j).
\label{eq:qm}
\ee

In summary, we have obtained: 

\begin{proposition}[Chebyshev method]
The Chebyshev method at level $m$ computes the polynomial $q_m$ of Eq.\ \eqref{eq:qm}, which provides the following upper and lower bounds to $d_p(M)$, 
\begin{align}
\label{eq:bound2}
\begin{split}
	&d_p(M)^p \geq \ntr(q_m(M)) - \Vert q_m-f_p  \Vert_\infty   \\
	 & d_p(M)^p \leq  \ntr(q_m(M)) +\Vert q_m-f_p \Vert_\infty.
\end{split}
\end{align}
\end{proposition}

\begin{figure}[htb]
\includegraphics[scale=0.6]{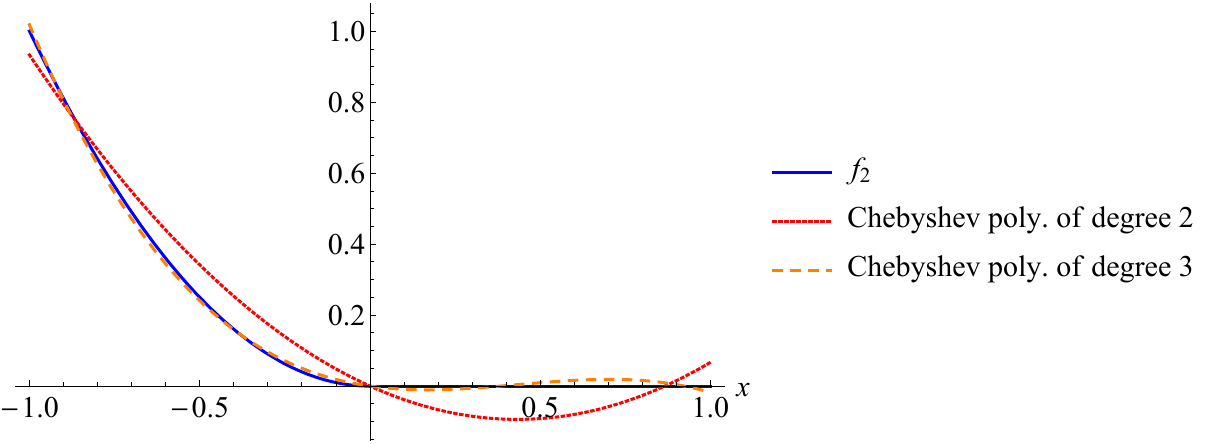}
\caption{The Chebyshev approximation of degree 2 and 3 of $f_2$. 
}\label{fig:app3}
\end{figure}

Note that by minimizing $\Vert q_m-f_p  \Vert_\infty$ one minimises $\Vert (q_m-f_p)_+  \Vert_\infty$ and $\Vert (q_m-f_p)_-  \Vert_\infty$ simultaneously. As an example, Figure \ref{fig:app3} shows two Chebyshev approximations of $f_2$.

\section{Numerical Results and Examples}
\label{sec:num}

We now discuss the numerical performance of the three methods presented in Section  \ref{sec:algo}. 
All computations were done with Mathematica on a standard computer with a 2.9GHz processor and 16 GB RAM. 
Our Mathematica code is available with this paper.

Since all three methods take the vector of the first $m$ normalized moments of a matrix $M$ as input, we first  study how long it takes to compute these moments for a  matrix of the form  \eqref{eq:Cl}. Specifically, we consider  random local matrices $A_j^{[i]}$ of size $2$, a sum length of $d=2$ and several values for the number $n$ of tensor factors. 
Note that this is a natural scenario in the tensor network paradigm, where one could imagine having a numerical tensor of the form  \eqref{eq:Cl} and having to decide, or give necessary or sufficient conditions, on its positive semidefiniteness.
 Table \ref{tab:mom} shows the average running time  to compute the first $m$ moments.
Note that already for $n=16$ we were not able to explicitly compute and store the matrix $M$, let alone decide its semidefiniteness on our computer.

\begin{table}[htb]\centering
\begin{tabular}{|c||*{5}{c|}}\hline
\diagbox{$n$}{$m$}
&\makebox[3em]{8}&\makebox[3em]{12}&\makebox[3em]{16}
\\\hline\hline
8 &0.02&0.60&13.82\\\hline
16 &0.09&1.44&28.01\\\hline
32 &0.18&2.60&55.52\\\hline
64 &0.33&4.88&111.44\\\hline
\end{tabular}
\caption{Average running time (in seconds) to compute the first $m$ moments for $n$ tensor factors for a matrix of the form of \eqref{eq:Cl}.}\label{tab:mom}
\end{table}

Now we study the running time of the three methods (with $p=2$) as a function of the number $m$ of moments used (Table \ref{tab:comp_avg}). Note that this does not include the computation of the moments, since they are part of the input. 
The moments have been produced from random instances of matrices of the form \eqref{eq:Cl} with local matrices $A_j^{[i]}$ of size $2$, $d=2$ and $n=8$. 
Note that the running time includes the full setup of all problems, as provided in the attached Mathematica file. 
In particular, for the Chebyshev method this includes the computation of the Chebyshev polynomial \eqref{eq:qm} of degrees $k=1,\ldots, m$ which interpolates $f_p$. These, however, do not depend on the input, and thus could be computed in advance. The algorithm itself  just has to compute the inner product of the coefficient vector and the moment vector, which  can be done in almost no time.

\begin{table}[htb]\centering
\begin{tabular}{|c||*{5}{c|}}\hline
\diagbox{$m$}{method}
&\makebox[7em]{sos polynomial}&\makebox[5em]{Handelman}&\makebox[5em]{Chebyshev}
\\\hline\hline
8   & 10.20 &0.16& 0.37\\\hline
12 &18.39 &0.40&0.76\\\hline
16 & 49.71&2.22&1.17\\\hline
32 &717.60 &23.85& 3.84\\\hline
\end{tabular}
\caption{Average running time (in seconds) of the three methods, using the first $m$ moments.}\label{tab:comp_avg}
\end{table}

To examine the qualitative performance of the three methods, we 
generated $10\,000$ random numbers uniformly in $[-\epsilon, 1]$, for several small values of $\epsilon>0$, and took them as eigenvalues of our matrix $M$. Note that the corresponding spectral measure is a close approximation of the uniform Lebesgue measure on this interval. 
We then computed the corresponding normalized moments, and checked how many moments each method needed to produce a positive lower bound on the distance to the psd cone, i.e.\ to detect non-positivity. 
Note that the smaller $\epsilon$, the harder the task.  
Table \ref{tab:qual} shows the average number of moments needed for each method, depending on $\epsilon$.
For $\epsilon=1/16$, the Chebyshev algorithm did never provide positive bounds when using  a number of moments that we could compute without running into numerical problems (for example when computing $\Vert f_p-q_m\Vert_\infty$, which is needed in the algorithm).

\begin{table}[htb]\centering
\begin{tabular}{|c||*{5}{c|}}\hline
\diagbox{$\epsilon$}{method}
&\makebox[7em]{sos polynomial}&\makebox[5em]{Handelman}&\makebox[5em]{Chebyshev}
\\\hline\hline
 1/2& 3&5&4\\\hline
 1/4&4 &7&8\\\hline
 1/8&5 &8&20\\\hline
 1/16& 6&17&?\\\hline
\end{tabular}
\caption{Average number of moments needed to detect non-positivity of a matrix with 10\, 000 random eigenvalues in the interval $[-\epsilon, 1]$.}\label{tab:qual}
\end{table}

Let us summarize and compare the three methods. Concerning the running time, the Chebyshev method is clearly the best. In addition, as mentioned above, its  running time can in practice be reduced to almost zero by computing the approximating polynomials beforehand. The Handelman method is also quite fast, in particular when compared to the sos polynomial method. On the other hand, the sos polynomial method needs significantly fewer moments  than the other methods in order to produce the same qualitative results. Computing many moments can also require a lot of time, depending on the format of the matrix. 

In order to compare both effects (running time versus number of moments needed), we conducted a final experiment. We produced random matrices again of the form  \eqref{eq:Cl} with local matrices of size $2$,  $d=2$ and different values for the number $n$ of tensor factors. For each method we first checked how many moments were needed to detect non-positivity, and then added the time to compute these moments to the actual running time of the method. Note that in practice one does not know in advance how many moments are needed to detect non-positivity. 
Interestingly,  the Chebyshev method falls behind the other two by far, due to the large number of moments it needs. 
The comparison of the Handelman and the sos polynomial method is summarized in  Table \ref{tab:comp_tot}. Note that also the Handelman method did not produce any meaningful result in the case of a very large matrix.

\begin{table}[htb]\centering
\begin{tabular}{|c||*{5}{c|}}\hline
\diagbox{$n$}{method}
&\makebox[7em]{sos polynomial }&\makebox[5em]{Handelman}
\\\hline\hline
16 & 0.63& 0.01\\\hline
24 & 1.33& 0.01\\\hline
32 & 7.34 & 1.46\\\hline
40 & 77 & ?\\\hline
\end{tabular}
\caption{Total running time (in seconds) to detect non-positivity, involving the time to compute the moments and the running time of the algorithm, for $n$ tensor factors.}\label{tab:comp_tot}
\end{table}

In summary, a general and clear decision  between the three methods based on their performance cannot be made. They differ greatly in terms of running time and in terms of the moments needed to obtain meaningful results.  However, the Chebyshev method seems to fall behind  the other two  in many relevant instances, at least for single matrices. If the matrices are not extremely large and the eigenvalues do not show an extreme behavior, the Handelman methods seems to perform best. The sos polynomial method can however  solve some of these extreme cases in which the other two methods fail.

\medskip

\emph{Files attached to this manuscript:} 
\begin{itemize}
\item\texttt{PSDBounds.nb}: This Mathematica package contains the three functions \texttt{sosApp}, \texttt{handelApp} and \texttt{chebyApp} that provide the upper and lower bounds, together with the approximating polynomials corresponding to the sos polynomial method, Handelman method, and Chebyshev method, respectively.
\item\texttt{TensorExample.nb}: This Mathematica notebook provides some examples to illustrate the use of the above package, in particular for matrices of the form \eqref{eq:Cl}.    
\end{itemize}

\emph{Acknowledgements.}---
	GDLC acknowledges funding of the Elise Richter Program of the FWF. TN acknowledges funding through the FWF project P 29496-N35 (free semialgebraic geometry and convexity). Most of this work was conducted while TF was at the Max Planck Institute for Mathematics in the Sciences. We thank Hilary Carteret and Andreas Thom for helpful comments on the topic.

\appendix

\section{Extension to von Neumann algebras}
\label{vNa}

\newcommand{\vNa}{\mathcal{N}}
\newcommand{\B}{\mathcal{B}}

Let $\vNa$ be a von Neumann algebra equipped with a faithful normal trace $\ntr$ satisfying $\ntr(1) = 1$. For example, $\vNa$ may be the group von Neumann algebra of a discrete group $\Gamma$, defined as the weak operator closure of the group algebra $\C[\Gamma]$ as acting on $\ell^2(\Gamma)$; the trace is given by $\ntr(x) := \langle e,xe\rangle$ with $e\in\Gamma$ being the unit. In the case where $\vNa = M_s(\C)$ is just the matrix algebra, the discussion presented here specializes to that of the main text.

For a Hermitian element $M\in\vNa$, the Schatten $p$-norm is again defined as
\[
	\| M \|_p := \left(\ntr(|M|^p)\right)^{1/p},
\]
where the absolute value and power functions are defined in terms of functional calculus. We can now ask the same question as in the main text: suppose we have $M$ of which we only know the values of the first $m$ moments
\[
	\ntr(M^k) \mbox{ for } k = 0,\ldots,m.
\]
Then what can we say about the $p$-distance from $M$ to the cone of positive elements?

It is straightforward to see that essentially all of the methods of the main text still apply without any change. The only differences are the following:

\begin{itemize}
\item Remark \ref{rem:scale} no longer applies, since the estimate that we used there involves the matrix size explicitly. Hence there is no bound on $\|M\|_\infty$ that could be computed from the moments. So in order to scale $M$ such that $\|M\|_\infty \leq 1$ is guaranteed, an a priori bound on $\|M\|_\infty$ needs to be known in addition to the moments. Otherwise our methods will not apply in their current form.
\item The proof of Proposition \ref{prop:psddist} is still essentially the same, but $\sigma$ needs to be generalized to the \emph{spectral scale}, and the corresponding inequality is~\cite[Corollary~3.3(1)]{Hi87}. \end{itemize}

Everything else is completely unchanged, including the fact that we are secretly addressing a version of the Hausdorff moment problem.


\end{document}